\documentclass{amsart}
\usepackage[utf8]{inputenc}
\usepackage{amssymb}
\usepackage{amsmath}
\usepackage{amsthm}
\usepackage{xcolor}
\usepackage{enumitem}
\usepackage{tikz}
\usepackage{hyperref}
\usepackage{mathtools}
\usepackage{comment}

\newcommand{\divi}{\textrm{div}}
\newcommand{\NN}{\mathbb{N}}
\newcommand{\Ccal}{\mathcal{C}}

\newtheorem{theorem}{Theorem}
\newtheorem*{theorem*}{Theorem}
\newtheorem{definition}[theorem]{Definition}
\newtheorem{remark}[theorem]{Remark}

\newtheorem{lemma}[theorem]{Lemma}
\newtheorem{proposition}[theorem]{Proposition}

\newtheorem{maintheorem}{Main Theorem}

\title{Weierstrass semigroups at ramified places of Frobenius nonclassical trinomial curves}
\author{Tiago Aprigio, João Paulo Guardieiro and Leonardo Oliveira}
\date{}
\address{Instituto de Ciências  Matemáticas e Computação, Av. Trabalhador São Carlense 400, 13.566-590 São Carlos-SP, Brazil.}
\email{tiagoaprigio@usp.br}
\address{Coordenação do Curso de Matemática - Bacharelado, Av. dos Portugueses, 1966, 65.080-805 São Luís-MA, Brazil.}
\email{joao.guardieiro@ufma.br}
\address{Instituto de Matemática, Estatística e Computação Científica da Universidade Estadual de
Campinas, R. Sergio Buarque de Holanda, 651, 13083-970 Campinas-SP, Brazil.}
\email{l237277@dac.unicamp.br}

\begin{document}

\begin{abstract}
   Frobenius nonclassical curves, introduced in 1985 by Stöhr and Voloch, are of significant interest due to their geometric properties and large number of rational points. In this paper we determine the Weierstrass semigroups at ramified places of the function fields of two families $\mathcal{C}_2$ and $\mathcal{C}_3$ of Frobenius nonclassical trinomial curves (that is, plane curves whose defining equations consist of exactly three monomials). For the first family, which is defined by a Kummer equation, we determine the semigroups at the places above the zero and pole of $x$ in the function field extension $\mathbb{F}_{q^n}(\mathcal{C}_2)/\mathbb{F}_{q^n}(x)$, as well as at a family of totally ramified places. For the second family, we determine the semigroups at all places of $\overline{\mathbb{F}}_{q^n}(\mathcal{C}_3)$ centered at the singular points of the plane model of the curves. For each of the two families, we present explicit functions whose poles are associated to the elements of the Weierstrass semigroups.
\end{abstract}


\maketitle

\section{Introduction}
During the storage and transmission of data, part of the information may be corrupted by noise in the channel. Error-correcting codes are tools used to retrieve the lost information, and, therefore, are very important for Information Theory. The interest of mathematicians in this topic may be traced back to the work of Hamming \cite{hamming}, in which the first linear code was constructed. A linear code over a finite field $\mathbb{F}_q$ with parameters $[n,k,d]$ is a $k$-dimensional
subspace of $\mathbb F_q^n$ whose nonzero elements have
at least $d$ nonzero coordinates, with equality for some element.

Roughly speaking, the dimension $k$ is related to the maximum number of words that can be transmitted with the code. The length $n$ is the size of each of these words; the larger $n$, the more difficult it is to implement the code in a computer, and the longer it takes to transmit the information. The ratio $k/n$ is called the rate of the code and measures, in a sense, its efficiency. A code with high rate carries more information per transmitted symbol. Finally, the minimum distance $d$ is related to the number of errors that the code can correct. It is then intuitive to think that a good error-correcting code has high minimum distance and small length. However, these three parameters are connected via the Singleton bound $d \leq n - k + 1$. Therefore, to construct codes with high minimum distance, it is also necessary to have a somewhat large length.

In 1981, Goppa considered the evaluation of functions on rational points of an algebraic curve over a finite field $\mathbb{F}_q$ to construct a linear code (see \cite{goppa}). His pioneering work showed that properties of the curve translate into information on the parameters of the associated code. This connection was the key step for the construction of a family of error-correcting codes based on algebraic curves whose parameters asymptotically exceed the Gilbert–Varshamov bound, see \cite{tsfasman} for more details, and this highlights the importance of a better understanding of curves over finite fields.

Goppa's construction uses Riemann–Roch spaces, and for one-point codes the parameters are governed by the Weierstrass semigroup at the chosen place. Given a rational place $P$ of the function field of a curve $\mathcal{X}$, the Weierstrass semigroup $H(P)$ at $P$ is the subset of $\mathbb{N} = \{0, 1, 2, \ldots\}$ composed of the numbers $n$ for which there exists a rational function $f$ whose pole divisor is precisely $(f)_{\infty} = nP$. Roughly speaking, a rational function can be written as the quotient of two polynomials, and we are interested in those functions whose denominator vanishes only at $P$ with ``multiplicity'' $n$. It can be shown (see \cite[Theorem 1.6.8]{stichtenoth}) that this set is closed under addition, and that its complement in the set of natural numbers, called the set of gaps at $P$, is composed of $g(\mathcal{X})$ elements, where $g(\mathcal{X})$ stands for the genus of $\mathcal{X}$. 

In 1986, Stöhr and Voloch studied algebraic curves using linear series on their function fields \cite{stohr-voloch}. Their work divided the curves into two categories: the Frobenius classical curves and the Frobenius nonclassical ones. Roughly speaking, Frobenius nonclassical curves tend to have more rational points than a typical curve of the same genus.  Codes constructed from them can therefore be long, which leaves room for a large minimum distance. 

In his PhD thesis \cite{guardieiro}, the second author of the present paper classified the curves that are given by a trinomial planar equation and that are Frobenius nonclassical with respect to the linear series induced by the degree-one Veronese morphism. He showed that, up to projective equivalence, there are only three families of such curves. 

In this work, we study the Weierstrass semigroups at some places of these curves. Determining the Weierstrass semigroup at every place of a function field is usually difficult, see for instance \cite{vicino1,vicino2,vicino3}. We therefore restrict ourselves to the ramified places in the field extension $K(x,y)/K(x)$, where $K$ is a suitable extension of $\mathbb{F}_q$ and $K(x,y)$ denotes the function field of these curves. 

Even with this restriction, we expect that the results of this paper can be useful for a future work in which we will construct error-correcting codes from these curves, and also for some classifying results of Frobenius nonclassical curves in terms of Weierstrass semigroups. Such results could follow from a comparison between the semigroups obtained here and those of the quadrinomial curves studied by the first and second authors in \cite{aprigio}.

The main results of this paper are the following. Let $q$ be a prime power and let $K=\mathbb F_{q^n}$. Assume first that $2\mid n$, set $t_2=(q^n-1)/(q^2-1)$, and consider the Kummer curve
\begin{equation}\label{tipo2}
 \mathcal C_2:\qquad y^{(q+1)t_2}=x^{t_2}(x^{(q-1)t_2}+1).
\end{equation}

\begin{maintheorem}
Let $P$ be any place of $K(\mathcal C_2)$ above the zero or the pole of $x$ in $K(x)$. 
Writing $L=q^2-1$, we have
$$
 H(P)=\big\langle (t_2-1)L+s :
                  s\in\langle q,q+1\rangle,\ 0\le s\le L\big\rangle .
$$
\end{maintheorem}

\begin{proof}
See Theorem \ref{semigroup}.
\end{proof}

\begin{maintheorem}
Let $\alpha\in K$ satisfy $\alpha^{(q-1)t_2}=-1$ and $P_\alpha$ be a place of $K(\mathcal C_2)$ above the zero of $x-\alpha$. We have
$$
 H(P_\alpha)=\big\langle (q+1)t_2-1,\ (q+1)t_2,\
                              (q^2-1)t_2-(q+1)\big\rangle .
$$
\end{maintheorem}

\begin{proof}
See Theorem \ref{semigroup-alpha}.
\end{proof}

Assume now that $3\mid n$, let $t_3=(q^n-1)/(q^3-1)$, $m=(q-1)t_3$, and $N=q^2+q+1$. Over $\overline K$, consider the curve
\begin{equation}\label{tipo3}
 \mathcal C_3:\qquad x^{(q+1)m}+x^m y^{(q+1)m}+y^m=0.
\end{equation}

\begin{maintheorem}
Let $\mathcal O$ be the set of places of $\overline K(\mathcal C_3)$ centered at $(1:0:0)$, $(0:1:0)$, or $(0:0:1)$ on the plane model. For every $P\in\mathcal O$,
$$
 H(P)=\big\langle (m-1)N+qn_1+(q+1)n_2 : \ n_1,n_2 \in\mathbb Z_{\ge0},\ n_1\le qn_2,\  qn_1+(q+1)n_2\le N\big\rangle .
$$

\end{maintheorem}

\begin{proof}
See Theorem \ref{semigroup3}.
\end{proof}

This paper is organized as follows: in Section \ref{sec:preliminaries} we establish some lemmas on numerical semigroups, and we provide the basic concepts we need throughout the text. In Section \ref{sec:C2} we study the curve $\mathcal{C}_2$ as in Equation \eqref{tipo2} and prove Main Theorems 1 and 2: we present some of its $\mathbb{F}_{q^n}$-automorphisms, the divisors of the coordinate functions and the Weierstrass semigroups at the ramified places of the function field of $\mathbb{F}_{q^n}(\mathcal{C}_2)$, by constructing explicit functions and using a genus count via Proposition \ref{lem:genus}. Finally, in Section \ref{sec:C3} we obtain analogous results for the curve $\mathcal{C}_3$ as in Equation \eqref{tipo3} and its ramified places of the function field $\overline{\mathbb{F}}_{q^n}(\mathcal{C}_3)$, and we prove Main Theorem 3.

\section{Preliminaries}\label{sec:preliminaries}
\subsection{Numerical semigroups}
As we mentioned before, given a rational place $P$ of the function field of a curve $\mathcal{X}$ and its Weierstrass semigroup $H(P)$, the sum of two elements of $H(P)$ is still an element of this set, and the cardinality of $\mathbb{N} \backslash H(P)$ is precisely the genus of the curve $\mathcal{X}$. Therefore, Weierstrass semigroups are important examples of numerical semigroups. For background on numerical semigroups and their applications, we refer the reader to \cite{numerical}. 

In this work, we construct the Weierstrass semigroups at some places by first obtaining rational functions with suitable pole divisors, and then showing that the genera of the associated numerical semigroups equal the genera of the curves. Proposition \ref{lem:genus} will be important for this final step, but first we fix the following notation.

\begin{itemize}

\item $[a,b]$ is the set of integers $n$ with $a\le n\le b$.
\item $\Lambda\subseteq\NN$ is a numerical semigroup, that is, $0\in\Lambda$, $\Lambda$ is closed under addition and $\NN\setminus\Lambda$ is finite. The elements of $\NN\setminus\Lambda$ are the gaps of $\Lambda$.
\item $g(\Lambda)$: genus of $\Lambda$, the number of gaps of $\Lambda$.
\item $\mu(\Lambda)$: multiplicity of $\Lambda$, the smallest positive element of $\Lambda$.
\item $c(\Lambda)$: conductor of $\Lambda$, the smallest integer $c$ such that every integer $n\ge c$ belongs to $\Lambda$.
\item $\langle T\rangle$: the semigroup generated by $T\subseteq\NN$, that is, the set of all finite sums of elements of $T$.
\item $\Sigma_k(T)\coloneqq\{s_1+\cdots+s_k \;:\; s_1,\ldots,s_k\in T\}$ for $k\ge1$, and $\Sigma_0(T)\coloneqq\{0\}$. Thus $\langle T\rangle=\bigcup_{k\ge0}\Sigma_k(T)$, and if $0\in T$, then $\Sigma_k(T)\subseteq\Sigma_{k+1}(T)$ for every $k\ge0$.
\end{itemize}

  \begin{lemma}\label{lem:Lambda_q}
Let $q\geq 2$ and consider the numerical semigroup
$$\Lambda_q := \langle\, q+1+q\cdot j \,:\, 0\leq j\leq q \,\rangle.$$
Then we have
\begin{enumerate}
    \item[(a)] $\Lambda_q=\{q\cdot n_1+(q+1)\cdot n_2 \,:\, n_1,n_2\in\mathbb{N},\ n_1\leq q\cdot n_2\}$.
    \item[(b)] The gaps of $\Lambda_q$ are the integers $r+jq$ with $1\leq r\leq q$ and $0\leq j\leq r-1$.
    \item[(c)] $\mu(\Lambda_q)=q+1$, $c(\Lambda_q)=q^2+1$ and $g(\Lambda_q)=\frac{q(q+1)}{2}$.
\end{enumerate}
\end{lemma}

\begin{proof}
Set $N:=q^2+q+1$ and, for $k\geq 0$, let $S_k$ be the set of sums of exactly $k$ generators (not necessarily distinct) of $\Lambda_q$, so that $S_0=\{0\}$ and $\Lambda_q=\bigcup_{k\geq 0}S_k$.

(a) An element of $S_k$ has the form
$$\sum_{\ell=1}^{k}(q+1+q\cdot j_\ell)=k(q+1)+q\,(j_1+\cdots+j_k), \qquad 0\leq j_1,\dots,j_k\leq q.$$
When $(j_1,\dots,j_k)$ runs over $[0,q]^k$, the sum $i:=j_1+\cdots+j_k$ takes every value in $[0,qk]$. Therefore
\begin{equation}\label{eq:Sk}
S_k=\{k(q+1)+qi \,:\, 0\leq i\leq qk\},
\end{equation}
that is, $S_k$ is the set of integers in $[k(q+1),kN]$ that are congruent to $k$ modulo $q$. This proves (a) with $n_1 = i$ and $n_2 = k$.

(b) Let $n\geq 1$ be an integer and let $r\in[1,q]$ be such that $n\equiv r\pmod q$. Then $n=r+jq$ for some $j\in\mathbb{N}$. We divide the proof into two cases.

\textbf{Case 1.} Suppose $0\leq j\leq r-1$. If $n\in\Lambda_q$, then $n\in S_k$ for some $k\geq 1$, and $k\equiv r\pmod q$ by Equation \eqref{eq:Sk}. Hence $k\geq r$ and
$$n\geq k(q+1)\geq r(q+1)>r+jq=n,$$
a contradiction. Thus $n$ is a gap of $\Lambda_q$.

\textbf{Case 2.} Suppose $j\geq r$, so that $n\geq r(q+1)$. Let $k\geq r$ with $k\equiv r\pmod q$. By Equation \eqref{eq:Sk}, $S_k$ contains every integer congruent to $r$ modulo $q$ in $[k(q+1),kN]$, and the first integer greater than $kN$ and congruent to $k$ modulo $q$ is $kN+q$. Since
$$(k+q)(q+1)=k(q+1)+q^2+q\leq k(q+1)+kq^2+q=kN+q,$$
the union of the sets $S_{r+\ell q}$, with $\ell\in\mathbb{N}$, contains every integer $m\geq r(q+1)$ with $m\equiv r\pmod q$. In particular, $n\in S_{r+\ell q}$ for some $\ell\in\mathbb{N}$, and so $n\in\Lambda_q$. This proves (b).

(c) From (b), for each $r\in[1,q]$ there are exactly $r$ gaps of $\Lambda_q$ congruent to $r$ modulo $q$. Hence
$$g(\Lambda_q)=\sum_{r=1}^{q} r=\frac{q(q+1)}{2}.$$
The largest gap congruent to $r$ modulo $q$ is $r+(r-1)q=r(q+1)-q$, which increases with $r$. So the largest gap of $\Lambda_q$ is $q(q+1)-q=q^2$, and therefore $c(\Lambda_q)=q^2+1$. Finally, the integers $1,\dots,q$ are gaps (take $j=0$ in (b)) and $q+1\in\Lambda_q$, so $\mu(\Lambda_q)=q+1$.
\end{proof}

\begin{lemma}\label{lem:sumset}
Let $\Lambda\neq\NN$ be a numerical semigroup, let $L$ be an integer with $L\ge c(\Lambda)+\mu(\Lambda)-1$, and set $T\coloneqq\Lambda\cap[0,L]$. Then
$$\Sigma_k(T)=\Lambda\cap[0,kL]\qquad\text{for every }k\ge1.$$
\end{lemma}
 
\begin{proof}
Write $c:=c(\Lambda)$ and $\mu:=\mu(\Lambda)$. Since $L\geq c$, both $0$ and $L$ belong to $T$. The inclusion $\Sigma_k(T)\subseteq\Lambda\cap[0,kL]$ is clear, since $\Lambda$ is closed under addition and every element of $T$ is at most $L$.

We prove the other inclusion by induction on $k$. For $k=1$, this is precisely the definition of $T$. Assume that $\Sigma_k(T)=\Lambda\cap[0,kL]$ for some $k\geq 1$, and let $n\in\Lambda\cap[0,(k+1)L]$. If $n\leq kL$, then $n\in\Sigma_k(T)\subseteq\Sigma_{k+1}(T)$, because $0\in T$. If $kL<n\leq (k+1)L$, then $n-(k-1)L\in[L,2L]$ and $(k-1)L\in\Sigma_{k-1}(T)$, because $L\in T$. Hence
$$n=(k-1)L+\big(n-(k-1)L\big)\in\Sigma_{k-1}(T)+[L,2L].$$
Since $\Sigma_{k-1}(T)+\Sigma_2(T)\subseteq\Sigma_{k+1}(T)$, it remains to show that $[L,2L]\subseteq\Sigma_2(T)$.

Let $m\in[L,2L]$ and write $m=L+j$, with $0\leq j\leq L$. If $j\in\Lambda$, then $j\in T$ and $m=L+j\in\Sigma_2(T)$. Otherwise, $j$ is a gap of $\Lambda$, so $j\leq c-1$. Let $s'$ be the smallest multiple of $\mu$ with $s'\geq j$. Then $s'\in\Lambda$ and
$$j\leq s'\leq j+\mu-1\leq c+\mu-2<L,$$
so $s'\in T$. Moreover, $s:=m-s'=L-(s'-j)$ satisfies $L-\mu+1\leq s\leq L$. Since $L-\mu+1\geq c$, we get $s\in\Lambda\cap[0,L]=T$. Therefore $m=s+s'\in\Sigma_2(T)$, which proves that $[L,2L]\subseteq\Sigma_2(T)$, and this concludes the proof.
\end{proof}

\begin{proposition}\label{lem:genus}
Let $\Lambda\neq\NN$ be a numerical semigroup, let $L$ be an integer with $L\ge c(\Lambda)+\mu(\Lambda)-1$, and set $T\coloneqq\Lambda\cap[0,L]$. For an integer $a\ge1$, let
$$H\coloneqq\big\langle (a-1)L+s \;:\; s\in T\big\rangle.$$
Then $H$ is a numerical semigroup with
$$g(H)=\frac{a(a-1)L}{2}-(a-1)+a\,g(\Lambda)\qquad\text{and}\qquad c(H)=a(a-1)L+c(\Lambda).$$
More precisely, if $a\ge2$ and $G$ denotes the set of gaps of $\Lambda$, then the set of gaps of $H$ is the disjoint union of the intervals
$$I_k\coloneqq\big[kaL+1,\ (k+1)(a-1)L-1\big],\qquad 0\le k\le a-2,$$
and of the sets $k(a-1)L+G$, for $1\le k\le a$.
\end{proposition}

\begin{proof}
The set $H$ is a numerical semigroup, since it contains $0$, it is closed under addition by definition, and its complement is finite, because it contains two consecutive integers (since $T$ also has this property).
Let $G$ be the set of gaps of $\Lambda$ and $c\coloneqq c(\Lambda)$. Recall that $G\subseteq[1,c-1]$, that $c-1\in G$, and that $L\ge c$.

\vspace{\baselineskip}
 
 \noindent \textbf{Case $a=1$.} Here $H=\langle T\rangle$. Since $T\subseteq\Lambda$, we have $H\subseteq\Lambda$. On the other hand, Lemma~\ref{lem:sumset} gives $\Lambda\cap[0,kL]=\Sigma_k(T)\subseteq H$ for every $k\ge1$, and so $\Lambda\subseteq H$. Hence $H=\Lambda$, and both formulas hold.
 
\vspace{\baselineskip}
 
 \noindent \textbf{Case $a\ge2$.} Let $m_0:=(a-1)L$. Since $0\in T$, $m_0$ is the smallest generator of $H$, and $m_0\geq L\geq c$. For $k\geq 1$, the sums of exactly $k$ generators of $H$ form the set
$$H_k:=\{(m_0+s_1)+\cdots+(m_0+s_k) : s_1,\dots,s_k\in T\}=km_0+\Sigma_k(T)=km_0+\big(\Lambda\cap[0,kL]\big),$$
where the last equality comes from Lemma~\ref{lem:sumset}. Setting $H_0:=\{0\}$, we have $H=\bigcup_{k\geq 0}H_k$. The set $H_k$ is contained in the interval
$$B_k:=\big[km_0,\,k(m_0+L)\big]=\big[k(a-1)L,\,kaL\big].$$
Since $G\subseteq[1,c-1]$ and $c\leq L\leq kL$, the integers in $[0,kL]$ that do not belong to $\Lambda$ are exactly the elements of $G$. Therefore, for $k\geq 1$,
\begin{equation}\label{eq:BkHk}
B_k\setminus H_k=km_0+G.
\end{equation}
By Equation \eqref{eq:BkHk} with $k=a$, the set $H_a\subseteq H$ contains all the integers in $[am_0+c,\,am_0+aL]$. This interval has $aL-c+1\geq m_0+1$ elements, because $L\geq c$. Since $m_0\in H$, adding multiples of $m_0$ to these integers we get every integer $n\geq am_0+c$. So every gap of $H$ is smaller than $am_0+c$.

Next we determine the gaps of $H$. For $j\geq a+1$ we have $H_j\subseteq B_j\subseteq[(a+1)m_0,\infty)$, and $(a+1)m_0\geq am_0+c$. Hence every element of $H$ smaller than $am_0+c$ lies in $B_0\cup\cdots\cup B_a$. Set
$$I_k:=\big[kaL+1,\,(k+1)(a-1)L-1\big],\qquad 0\leq k\leq a-2.$$
The intervals
$$B_0=\{0\},\ I_0,\ B_1,\ I_1,\ \dots,\ B_{a-2},\ I_{a-2},\ B_{a-1},\ B_a$$
are consecutive and pairwise disjoint, except for $B_{a-1}$ and $B_a$, which share the endpoint $am_0\in H$, and their union is precisely $[0,a^2L]$. Therefore the elements of the intervals $I_k$ are gaps of $H$. Moreover, for $1\leq k\leq a$, the set $km_0+G$ lies in $B_k\setminus\{am_0\}$, because $G\subseteq[1,c-1]$ and $m_0\geq c$, so its elements are gaps of $H$ by Equation \eqref{eq:BkHk}, while every element of $B_k\setminus(km_0+G)$ belongs to $H_k\subseteq H$. Hence the gaps of $H$ are exactly the elements of the pairwise disjoint sets $I_k$, with $0\leq k\leq a-2$, and $km_0+G$, with $1\leq k\leq a$. Since each $I_k$ has $(a-1-k)L-1$ elements, we get
$$g(H)=\sum_{k=0}^{a-2}\big((a-1-k)L-1\big)+a|G|=L\sum_{j=1}^{a-1}j-(a-1)+a\,g(\Lambda)=\frac{a(a-1)L}{2}-(a-1)+a\,g(\Lambda).$$
Finally, the largest gap of $H$ is $am_0+(c-1)$, because $c-1\in G$ and every integer $n\geq am_0+c$ belongs to $H$. Hence $c(H)=am_0+c=a(a-1)L+c(\Lambda)$.
\end{proof}

\subsection{Ramified places of Frobenius nonclassical trinomial curves} Let $p$ be a prime number and $q = p^e$, with $e \geq 1$. Frobenius nonclassical curves (with respect to the linear series induced by the degree-one Veronese morphism) are irreducible curves over a finite field $\mathbb{F}_{q^n}$ with the following geometric property: the image of any smooth point under the $\mathbb{F}_{q^n}$-Frobenius map lies in its tangent line (for a proper definition of such curves, we refer the reader to \cite{stohr-voloch}). 

In his PhD thesis \cite{guardieiro}, the second author of the present paper showed that, up to $\mathbb{F}_{q^n}$-projective equivalence, there are only three families of Frobenius nonclassical curves whose plane models are given by a trinomial equation:
$$\mathcal{C}_1: x^{t_1} + y^{t_1} + 1 = 0, \qquad {t_1} = \frac{q^n - 1}{q - 1}
$$
$$\mathcal{C}_2: y^{(q+1){t_2}} = x^{q{t_2}} + x^{t_2}, \qquad {t_2} = \frac{q^n - 1}{q^2 - 1}, \qquad 2 \mid n
$$
$$\mathcal{C}_3: x^{(q^2-1){t_3}} + x^{(q-1){t_3}}y^{(q^2-1){t_3}} + y^{(q-1){t_3}} = 0, \qquad {t_3} = \frac{q^n - 1}{q^3 - 1}, \qquad 3 \mid n
$$
Observe that $t_i$ is an integer, since $i \mid n$, and that $t_i \equiv 1 \pmod{p}$, for $i \in \{1,2, 3\}$. 

\begin{remark}
    It follows from \cite[Propositions 2.84, 2.85 and 2.86]{guardieiro} that the curve $\mathcal{C}_1$ is smooth, that $\mathcal{C}_2$ is smooth if, and only if, $t_2 = 1$ (in which case the curve is precisely the Hermitian curve over $\mathbb{F}_{q^2}$) and that $\mathcal{C}_3$ is smooth if, and only if, $(q,n) = (2,3)$ (in which case the curve is precisely the Klein quartic over $\mathbb{F}_8$). The Weierstrass semigroups at the rational points of such smooth curves $\mathcal{C}_2$ and $\mathcal{C}_3$ are well-known, and coincide with the ones obtained here.
\end{remark}

The genus and number of rational places of these curves were also obtained in the thesis (see Theorem B of Chapter 3).

\begin{proposition}\label{curves-places}
    Consider the curves $\mathcal{C}_1$, $\mathcal{C}_2$ and $\mathcal{C}_3$ as above. We have the following.
    \begin{enumerate}
        \item The curve $\mathcal{C}_1$ has genus $g(\mathcal{C}_1) = \frac{t_1^2 - 3t_1 + 2}{2}$, and the number of $\mathbb{F}_{q^n}$-rational places of its function field equals $t_1^2(q - 2) + 3t_1$.
        \item The curve $\mathcal{C}_2$ has genus $g(\mathcal{C}_2) = \frac{(q^2 - 1)t_2^2 - (q + 1)t_2 + 2}{2}$, and the number of $\mathbb{F}_{q^n}$-rational places of its function field equals $t_2^2(q^3 - q) + t_2(q + 1)$.
        \item The curve $\mathcal{C}_3$ has genus $g(\mathcal{C}_3) = \frac{(q^3 - 1)(q - 1)t_3^2 - 3(q - 1)t_3 + 2}{2}$, and the number of $\mathbb{F}_{q^n}$-rational places of its function field equals 
        $$\left\{\begin{array}{rr}
            t_3^2(q^2 - 1)(q^3 - 1) + 3t_3(q - 1) & \textrm{if } p = 2 \\
            t_3^2(q^2 - 1)(q^3 - 1) & \textrm{if } p \neq 2
        \end{array}\right..$$
    \end{enumerate}
\end{proposition}

The goal of this paper is to determine the Weierstrass semigroups at the ramified places in the field extensions $\mathbb{F}_{q^n}(\mathcal{C}_i)/\mathbb{F}_{q^n}(x)$, for $i \in \{1, 2, 3\}$. The first step is then to describe precisely such places. Note that the function fields of the curves $\mathcal{C}_1$ and $\mathcal{C}_2$ are Kummer extensions of $\mathbb{F}_{q^n}(x)$, since their planar model can be written as $y^{\ell} = f(x)$, for suitable $\ell \in \mathbb{N}$ not divisible by the characteristic $p$ of the field and $f(x) \in \mathbb{F}_{q^n}[x]$ with at least one simple root. From \cite[Proposition 3.7.3]{stichtenoth}, we then conclude that the only (possibly) ramified places at the field extension are those associated to the pole of $x$ and to the roots of $f(x)$.

Since the polynomial $x^{t_1} + 1$ is of degree $t_1 = [\mathbb{F}_{q^n}(\mathcal{C}_1) : \mathbb{F}_{q^n}(x)]$ and is separable in $\mathbb{F}_{q^n}$, the only ramified places at the field extension $\mathbb{F}_{q^n}(\mathcal{C}_1)/\mathbb{F}_{q^n}(x)$ are those associated to the roots of this polynomial. Furthermore, these places are totally ramified and $\mathcal{C}_1$ belongs to the family of Fermat curves. A result of Hasse \cite{hasse} now implies that the Weierstrass semigroups at each of these places equals $\langle t_1 - 1, t_1 \rangle$, and so we do not consider the curve $\mathcal{C}_1$ in this work. We should mention that the Weierstrass semigroups at the remaining places of the function field of a Fermat curve are also relevant, and we are already studying them (in the case of a Frobenius nonclassical Fermat curve).

For the curve $\mathcal{C}_2$, we conclude that the zero and the pole of $x$ each split into $t_2$ places of $\mathbb F_{q^n}(\mathcal C_2)$, each with ramification index $q+1$. Also, the zero of $x - \alpha \in \mathbb{F}_{q^n}(x)$, with $\alpha^{(q - 1)t_2} + 1 = 0$, totally ramifies into a single place $P_{\alpha}$ of $\mathbb{F}_{q^n}(\mathcal{C}_2)$, with ramification index $(q + 1)t_2$. The Fundamental Equality \cite[Theorem 3.1.11]{stichtenoth} shows that all these places are rational. Every other rational place of $\mathbb{F}_{q^n}(x)$ splits completely in the field extension $\mathbb{F}_{q^n}(\mathcal{C}_2)/\mathbb{F}_{q^n}(x)$. The classification of Weierstrass semigroups in Kummer extensions is a well-studied topic, but, in opposition to our paper, it is common to fully determine it only for totally ramified places (see, for instance, \cite{abq2019,castellanos2024,yang-hu,zhang-zhao}). We also obtain the semigroups by presenting explicit functions that give rise to the non-gap numbers, instead of relying on theoretical results.

Finally, it was shown in \cite[Theorem 3.20]{guardieiro} that there are $(q - 1)t_3$ places of $\overline{\mathbb{F}}_{q^n}(\mathcal{C}_3)$ lying over the pole of $x$, and $2(q - 1)t_3$ places lying over its zero. Furthermore, these places are $\mathbb{F}_{q^{2n}}$-rational, and they are $\mathbb{F}_{q^n}$-rational if and only if $p = 2$. Therefore, when studying the curve $\mathcal{C}_3$, we consider the function field $\overline{\mathbb{F}}_{q^n}(\mathcal{C}_3)$, of functions with coefficients in the algebraic closure of $\mathbb{F}_{q^n}$, so we do not need to worry about the degree of the places. Any other place of $\overline{\mathbb{F}}_{q^n}(x)$ splits completely in the function field extension.

\begin{proposition}\label{non-singular-points}
    Let $a\in \overline{\mathbb{F}}_{q^n}^*$, and denote by $P_a$ the unique zero of $x-a$ on the function field $\overline{\mathbb{F}}_{q^n}(x)$. Then, $P_a$ splits completely in the field extension $\overline{\mathbb{F}}_{q^n}(\mathcal{C}_3) / \overline{\mathbb{F}}_{q^n}(x)$.
\end{proposition}
\begin{proof}
    Consider the polynomial
    $$W(T) = a^{(q^2 - 1)t_3} + a^{(q - 1)t_3}T^{(q^2 - 1)t_3} + T^{(q - 1)t_3}.$$
    An analysis of the derivative of $W$ shows that $W$ has only simple roots in $\overline{\mathbb{F}}_{q^n}$, which we denote by $b_1, \ldots, b_{(q^2 - 1)t_3}$. We also define $W_i(T) = T - b_i$, for $i = 1, \ldots, (q^2 - 1)t_3$.

    Since $v_{P_a}(x) = 0$, the defining equation of $\mathcal{C}_3$ can be seen as the minimal polynomial of $y$ over $\overline{\mathbb{F}}_{q^n}(x)$, with coefficients in $\mathcal{O}_{P_a}$. The decomposition of this polynomial over the residue field of $\mathcal{O}_{P_a}$ is given by
    $$\prod_{i = 1}^{(q^2 - 1)t_3}(T - b_i) = \prod_{i = 1}^{(q^2 - 1)t_3}\overline{W}_i(T).$$
    From \cite[Theorem 3.3.7]{stichtenoth}, we then conclude that there are $(q^2 - 1)t_3$ places of $\overline{\mathbb{F}}_{q^n}(x,y)$ lying over $P_a$, each one given by the common zero of $x - a$ and $y - b_i$. The ramification index of each one of these places is then 1.
\end{proof}

Therefore, each point $(a : b : 1) \in \mathcal{C}_3$, with $a \neq 0$ (and, therefore, $b \neq 0$) is associated to a single place $P_{(a,b)}$ of $\overline{\mathbb{F}}_{q^n}(\mathcal{C}_3)$. 

\section{The curve \texorpdfstring{$\mathcal{C}_2$}{C2}}\label{sec:C2}
We already mentioned that the zero of $x$ in $\mathbb{F}_{q^n}(x)$ corresponds to $t_2$ places in the function field $\mathbb{F}_{q^n}(\mathcal{C}_2)$, which we denote by $P_0^1, \ldots, P_0^{t_2}$. The ramification index of each place is $q + 1$. Analogously, the $t_2$ places of $\mathbb{F}_{q^n}(\mathcal{C}_2)$ associated to the pole of $x$ are denoted by $P_{\infty}^1, \ldots, P_{\infty}^{t_2}$, and their ramification index is also $q + 1$. We denote by $P_{\alpha}$, with $\alpha$ a root of $x^{(q - 1)t_2} + 1$, the only place in $\mathbb{F}_{q^n}(\mathcal{C}_2)$ associated to the zero of $x - \alpha$. Finally, any other rational place of $\mathbb{F}_{q^n}(\mathcal{C}_2)$ can be associated univocally to a rational point in the curve. We denote by $P_{(a,b)}$ the place associated to $(a : b : 1) \in \mathcal{C}_2$.

\subsection{The places \texorpdfstring{$P_0^i$}{P0i} and \texorpdfstring{$P_{\infty}^j$}{Pj}}
We first show that the Weierstrass semigroups at the places $P_0^i$ and $P_{\infty}^i$, for $i \in \{1, \ldots, t_2\}$ are the same. Then we find the semigroup at $P_0^1$ by presenting explicit functions giving rise to the pole numbers.

\begin{proposition}\label{Weierstrass-semigroups}
    The Weierstrass semigroups at the places $P_0^1, \ldots, P_0^{t_2}, P_{\infty}^1, \ldots, P_{\infty}^{t_2}$ are the same.
\end{proposition}
\begin{proof}
    Observe that the maps
    $$\theta_{\gamma,\delta}(x,y) = (\gamma x,\delta y), \qquad \gamma^{t_2} = \delta^{(q+1)t_2} \in \mathbb{F}_q$$
    and
    $$\theta_2(x,y) = \left(\frac{1}{x}, \frac{y}{x}\right)$$
    are automorphisms of $\mathcal{C}_2$.

    In particular, since $(q + 1)t_2 = \frac{q^n - 1}{q - 1}$ is a divisor of $q^n - 1$, the field $\mathbb{F}_{q^n}$ contains all $(q  + 1)t_2$-th roots of unity. The automorphisms $\theta_{1, \delta}$, with $\delta^{(q + 1)t_2} = 1$, are then the elements of the Galois group of the field extension $\mathbb{F}_{q^n}(\mathcal{C}_2) / \mathbb{F}_{q^n}(x)$ (this is known from Kummer theory).

    It is known that the Galois group of the field extension $\mathbb{F}_{q^n}(\mathcal{C}_2) / \mathbb{F}_{q^n}(x)$ acts transitively on the extensions of a place of $\mathbb{F}_{q^n}(x)$ (see \cite[Theorem 3.7.1]{stichtenoth}). Therefore, we conclude that
    $$\{P_0^1, \ldots, P_0^{t_2}\} \quad \textrm{and} \quad \{P_{\infty}^1, \ldots, P_{\infty}^{t_2}\}$$
    are orbits under the action of the Galois group of $\mathbb{F}_{q^n}(\mathcal{C}_2)/\mathbb{F}_{q^n}(x)$.
    
    The automorphism $\theta_2$ maps the point $(0:0:1)$ to $(1:0:0)$, and so maps a place $P_0^i$ to a place $P_{\infty}^j$. We then conclude that
    $$\{P_0^1, \ldots, P_0^{t_2}, P_{\infty}^1, \ldots, P_{\infty}^{t_2}\}$$
    is contained in an orbit under the action of the automorphism group of $\mathcal{C}_2$.

    It follows from \cite[Lemma 3.5.2(a)]{stichtenoth} that if two places lie on the same orbit under the action of the automorphism group, they have the same Weierstrass semigroup. Therefore, we conclude that $$H(P_0^i)=H(P_{\infty}^j), \quad \textrm{for every } i,j\in\{1,\ldots,t_2\}.$$
\end{proof}

To construct explicit functions that give us the Weierstrass semigroup at $P_0^1$, we first need to find the divisors of the coordinate functions $x$ and $y$.

\begin{proposition}\label{divisors-x-y}
    Let $\alpha$ be a root of $x^{(q - 1)t_2} + 1 = 0$. The principal divisors of the functions $x$, $x - \alpha$ and $y$ are given by
    $$\divi(x) = (q+1)\sum_{i=1}^{t_2}P_0^i - (q+1)\sum_{i=1}^{t_2}P_{\infty}^i$$
    $$\divi(x - \alpha) = (q + 1)t_2P_{\alpha} - (q+1)\sum_{i=1}^{t_2}P_{\infty}^i$$
    and
    $$\divi(y) = \sum_{i=1}^{t_2}P_0^i + \sum_{\alpha^{(q-1)t_2}+1=0}P_{\alpha} - q\sum_{i=1}^{t_2}P_{\infty}^i$$
\end{proposition}
\begin{proof}
    The expressions for $\divi(x)$ and $\divi(x - \alpha)$ come from the definition of the ramification index and the discussion at the beginning of this section. To obtain the expression for $\divi(y)$, we note that, from the equation of the curve $\mathcal{C}_2$,
    $$y^{(q+1)t_2} = x^{t_2}(x^{(q-1)t_2} + 1) = x^{t_2} \cdot \prod_{\alpha^{(q-1)t_2} + 1 = 0}(x - \alpha).$$ 
    Therefore,
    $$(q+1)t_2 \cdot \divi(y) = t_2 \cdot \divi(x) + \sum_{\alpha^{(q-1)t_2} + 1 = 0}\divi(x-\alpha),$$
    and this yields the desired result.
\end{proof}

\begin{lemma}\label{lem:w}
Let $w:=y^{q+1}/x\in\mathbb{F}_{q^n}(\mathcal{C}_2)$ and, for $i=1,\dots,t_2$, let $\beta_i:=w(P_0^i)$ be the residue class of $w$ at $P_0^i$. Then the following hold.
\begin{enumerate}[label=\upshape(\roman*)]
\item The divisor of $w$ is
$$\operatorname{div}(w)=(q+1)\sum_{\alpha^{(q-1)t_2}+1=0}P_\alpha-(q^2-1)\sum_{j=1}^{t_2}P_\infty^j.$$
\item The elements $\beta_1,\dots,\beta_{t_2}$ are pairwise distinct, and they are exactly the $t_2$-th roots of unity in $\mathbb{F}_{q^n}$.
\item For every $i\in\{1,\dots,t_2\}$ we have
$$\operatorname{div}(w-\beta_i)=(q^2-1)t_2\,P_0^i-(q^2-1)\sum_{j=1}^{t_2}P_\infty^j.$$
\end{enumerate}
\end{lemma}

\begin{proof}
By Proposition~\ref{divisors-x-y} we have
$$\operatorname{div}(w)=(q+1)\operatorname{div}(y)-\operatorname{div}(x)=(q+1)\sum_{\alpha^{(q-1)t_2}+1=0}P_\alpha-(q^2-1)\sum_{j=1}^{t_2}P_\infty^j,$$
which proves (i). In particular $w$ is a unit at each $P_0^i$, and the pole divisor of $w$ is $(q^2-1)\sum_{j=1}^{t_2}P_\infty^j$, which has degree $(q^2-1)t_2$.

The defining equation of $\mathcal{C}_2$ gives
$$w^{t_2}=\frac{y^{(q+1)t_2}}{x^{t_2}}=\frac{x^{qt_2}+x^{t_2}}{x^{t_2}}=x^{(q-1)t_2}+1,$$
and therefore, again by Proposition~\ref{divisors-x-y},
$$\operatorname{div}\big(w^{t_2}-1\big)=(q-1)t_2\operatorname{div}(x)=(q^2-1)t_2\sum_{i=1}^{t_2}P_0^i-(q^2-1)t_2\sum_{j=1}^{t_2}P_\infty^j.$$

Let $\mu_{t_2}$ be the set of $t_2$-th roots of unity in $\overline{\mathbb{F}}_{q^n}$. Since $t_2$ divides $q^n-1$, the set $\mu_{t_2}$ is contained in $\mathbb{F}_{q^n}$ and has exactly $t_2$ elements, so that
$$w^{t_2}-1=\prod_{\zeta\in\mu_{t_2}}(w-\zeta).$$
For each $\zeta\in\mu_{t_2}$ the function $w-\zeta$ has the same pole divisor as $w$, hence its zero divisor has degree $(q^2-1)t_2$. Every zero of $w-\zeta$ is a zero of $w^{t_2}-1$, and thus belongs to $\{P_0^1,\dots,P_0^{t_2}\}$. Moreover, for $\zeta\neq\zeta'$ the functions $w-\zeta$ and $w-\zeta'$ have no common zero, since their difference is the nonzero constant $\zeta'-\zeta$.

Fix $i\in\{1,\dots,t_2\}$. Since $v_{P_0^i}(w^{t_2}-1)>0$, exactly one factor $w-\zeta$ vanishes at $P_0^i$, and all the other factors are units at $P_0^i$. This factor is necessarily $w-\beta_i$, and
$$v_{P_0^i}(w-\beta_i)=v_{P_0^i}\big(w^{t_2}-1\big)=(q^2-1)t_2.$$
Consequently, for each $\zeta\in\mu_{t_2}$ the zero divisor of $w-\zeta$ equals $(q^2-1)t_2\sum_{\beta_i=\zeta}P_0^i$. Comparing degrees, we conclude that there is exactly one index $i$ with $\beta_i=\zeta$. Hence $i\mapsto\beta_i$ is a bijection from $\{1,\dots,t_2\}$ onto $\mu_{t_2}$, which proves (ii), and the zero divisor of $w-\beta_i$ is $(q^2-1)t_2\,P_0^i$, which proves (iii).
\end{proof}

\begin{proposition}\label{prop:f}
Let $w=y^{q+1}/x$ and $\beta:=\beta_1$ be as in Lemma~\ref{lem:w}, and define
$$f:=\frac{1}{x^{(q-1)(t_2-1)}}\left[1+\beta^{-1}w+\big(\beta^{-1}w\big)^2+\cdots+\big(\beta^{-1}w\big)^{t_2-1}\right].$$
Then
$$\operatorname{div}(f)=(q^2-1)\sum_{i=2}^{t_2}P_0^i-(q^2-1)(t_2-1)P_0^1.$$
\end{proposition}

\begin{proof}
By Proposition~\ref{divisors-x-y} we have
$$\operatorname{div}\left(\frac{1}{x^{(q-1)(t_2-1)}}\right)=(q^2-1)(t_2-1)\sum_{j=1}^{t_2}P_\infty^j-(q^2-1)(t_2-1)\sum_{i=1}^{t_2}P_0^i.$$
Let $g:=1+\beta^{-1}w+\cdots+\big(\beta^{-1}w\big)^{t_2-1}$, so that $f=g/x^{(q-1)(t_2-1)}$ and
\begin{equation}\label{eq:factor}
\big(\beta^{-1}w\big)^{t_2}-1=\big(\beta^{-1}w-1\big)\,g.
\end{equation}
Since $\beta^{t_2}=1$ by Lemma~\ref{lem:w}(ii), the defining equation of $\mathcal{C}_2$ gives
$$\big(\beta^{-1}w\big)^{t_2}-1=\frac{y^{(q+1)t_2}}{x^{t_2}}-1=\frac{x^{qt_2}+x^{t_2}}{x^{t_2}}-1=x^{(q-1)t_2},$$
and hence, by Proposition~\ref{divisors-x-y},
$$\operatorname{div}\Big(\big(\beta^{-1}w\big)^{t_2}-1\Big)=(q^2-1)t_2\sum_{i=1}^{t_2}P_0^i-(q^2-1)t_2\sum_{j=1}^{t_2}P_\infty^j.$$
On the other hand, $\beta^{-1}w-1=\beta^{-1}(w-\beta_1)$, so Lemma~\ref{lem:w}(iii) yields
$$\operatorname{div}\big(\beta^{-1}w-1\big)=(q^2-1)t_2\,P_0^1-(q^2-1)\sum_{j=1}^{t_2}P_\infty^j.$$
It follows from \eqref{eq:factor} that
$$\operatorname{div}(g)=(q^2-1)t_2\sum_{i=2}^{t_2}P_0^i-(q^2-1)(t_2-1)\sum_{j=1}^{t_2}P_\infty^j.$$
Adding the divisor of $1/x^{(q-1)(t_2-1)}$ computed above, we obtain
$$\operatorname{div}(f)=(q^2-1)\sum_{i=2}^{t_2}P_0^i-(q^2-1)(t_2-1)P_0^1,$$
which completes the proof.
\end{proof}

\begin{theorem}\label{semigroup}
    For any $i\in\{1,\ldots,t_2\}$, we have
    $$H(P_0^i)=H(P_\infty^i)=\big\langle (q^2-1)(t_2-1)+qn_1+(q+1)n_2 \;:\; n_1,n_2\in\NN,\ qn_1+(q+1)n_2\le q^2-1\big\rangle.$$
\end{theorem}
\begin{proof}
By Proposition \ref{Weierstrass-semigroups}, it is enough to consider the place $P_0^1$. Denote by $H$ the semigroup on the right-hand side.
 
Let $n_1,n_2\in\NN$ with $qn_1+(q+1)n_2\le q^2-1$, and let $f$ be as in Proposition \ref{prop:f}. From Propositions \ref{divisors-x-y} and \ref{prop:f}, the divisor of $\frac{y^{n_1}}{x^{n_1}}\cdot\frac{1}{x^{n_2}}\cdot f$ is
$$\begin{aligned}
&\big[(q^2-1)-qn_1-(q+1)n_2\big]\sum_{i=2}^{t_2}P_0^i+n_1\sum_{\alpha^{(q-1)t_2}+1=0}P_\alpha+\big[n_1+(q+1)n_2\big]\sum_{i=1}^{t_2}P_\infty^i\\
&\qquad-\big[(q^2-1)(t_2-1)+qn_1+(q+1)n_2\big]P_0^1.
\end{aligned}$$
Since the first three coefficients are nonnegative, $P_0^1$ is the only pole of this function, and its order is $(q^2-1)(t_2-1)+qn_1+(q+1)n_2$. Hence $H\subseteq H(P_0^1)$.
 
Now let $\Lambda\coloneqq\langle q,q+1\rangle$ and $L\coloneqq q^2-1$. It is known (see, for instance, \cite{sylvester1884}) that $c(\Lambda)+\mu(\Lambda)-1=(q^2-q)+q-1=L$. Moreover,
$$T\coloneqq\Lambda\cap[0,L]=\{qn_1+(q+1)n_2 \;:\; n_1,n_2\in\NN,\ qn_1+(q+1)n_2\le q^2-1\},$$
and so $H=\langle (t_2-1)L+s \;:\; s\in T\rangle$. Applying Proposition~\ref{lem:genus} with $a=t_2$, we get
$$g(H)=\frac{t_2(t_2-1)(q^2-1)}{2}-(t_2-1)+t_2\cdot\frac{q^2-q}{2}=\frac{(q^2-1)t_2^2-(q+1)t_2+2}{2}=g(\Ccal_2),$$
where the last equality is stated in Proposition \ref{curves-places}. Since $H\subseteq H(P_0^1)$ and both semigroups have genus $g(\Ccal_2)$ we have $H=H(P_0^1)$.
\end{proof} 

\subsection{The points \texorpdfstring{$P_{\alpha}$}{Pa}}
For each root $\gamma$ of $x^{(q-1)t_2}=1$ (in particular, $\gamma^{t_2} \in \mathbb{F}_q^*$), the surjectivity of the norm function ensures that there exists $\delta \in \mathbb{F}_{q^n}^*$ such that $\delta^{(q + 1)t_2} = \gamma^{t_2}$. The map $(x, y) \mapsto (\gamma x, \delta y)$ is now an automorphism of $\mathcal{C}_2$, and it maps $P_\alpha$ to $P_{\gamma\alpha}$. Hence all places $P_\alpha$ have the same Weierstrass semigroup, which we now determine. The semigroup at those places belongs to the class of telescopic semigroups.

\begin{definition}\cite[Definition 6.1]{hlp1998}
    Let $(a_1, \ldots, a_k)$ be a sequence of natural numbers such that $\gcd(a_1, \ldots, a_k) = 1$. For $i = 1, \ldots, k$, define
    $$d_i \coloneqq \gcd(a_1, \ldots, a_i), \quad A_i \coloneqq \left\{\frac{a_1}{d_i}, \ldots, \frac{a_i}{d_i}\right\}$$ and $H_i$ the numerical semigroup generated by $A_i$. The sequence $(a_1, \ldots, a_k)$ is called a \textbf{telescopic sequence} if $\frac{a_i}{d_i} \in H_{i-1}$ for every $i = 2, \ldots, k$. The numerical semigroup $H = \langle a_1, \ldots, a_k \rangle$ is called a \textbf{telescopic semigroup} if $(a_1, \ldots, a_k)$ is a telescopic sequence.
\end{definition}

\begin{proposition}\cite[Lemma 6.5]{hlp1998}
    Let $H = \langle a_1, \ldots, a_k \rangle$ be a telescopic semigroup. Then, the genus of $H$ is given by
    $$g(H) = \frac{1}{2}\left(1 - a_1 + \sum_{i = 2}^k\left(\frac{d_{i-1}}{d_i} - 1\right)a_i \right).$$
\end{proposition}

\begin{theorem}\label{semigroup-alpha}
    Let $\alpha \in \mathbb{F}_{q^n}$ be any root of $x^{(q - 1)t_2} + 1$. Then $$H(P_{\alpha}) = \langle (q + 1)t_2 - 1, (q + 1)t_2, (q^2 - 1)t_2 - (q + 1)\rangle.$$
\end{theorem}
\begin{proof}
    We use Proposition \ref{divisors-x-y} to conclude that
    $$\operatorname{div}\Big(\frac{y}{x-\alpha}\Big)=\sum_{i = 1}^{t_2} P_0^i+\sum_{\substack{\alpha'^{(q - 1)t_2} + 1 = 0 \\\alpha'\neq\alpha}}P_{\alpha'}+\sum_{i = 1}^{t_2} P_\infty^i-\big((q+1)t_2 - 1\big)P_\alpha$$
    $$\operatorname{div}\Big(\frac{1}{x-\alpha}\Big)=(q+1)\sum_{i = 1}^{t_2} P_\infty^i - (q+1)t_2\,P_\alpha$$
    and
    $$\operatorname{div}\Big(\frac{y^{q+1}}{x(x-\alpha)^{q-1}}\Big)= (q+1)\sum_{\substack{\alpha'^{(q - 1)t_2} + 1 = 0 \\\alpha'\neq\alpha}}P_{\alpha'}-\big((q^2-1)t_2-(q+1)\big)P_\alpha$$
    Therefore, we conclude that
    $$H\coloneqq\langle (q + 1)t_2 - 1, (q + 1)t_2, (q^2 - 1)t_2 - (q + 1) \rangle \subseteq H(P_{\alpha}).$$ 
    To prove equality, we just need to show that $g(H) = g(\mathcal{C}_2)$.
    
    Defining $a_1 = (q + 1)t_2, a_2 = (q^2 - 1)t_2 - (q + 1)$ and $a_3 = (q + 1)t_2 - 1$, we can see that $H$ is a telescopic semigroup generated by $a_1,a_2,a_3$. Indeed, following the previous notation, we have
    $$d_1 = (q + 1)t_2, \quad d_2 = q + 1, \quad d_3 = 1, \quad H_1 = \langle 1\rangle, \quad H_2 = \langle t_2, (q - 1)t_2 - 1\rangle.$$
    Since
    $$\frac{a_2}{d_2} = (q - 1)t_2 - 1 \in H_1 \quad \textrm{and} \quad \frac{a_3}{d_3} = (q + 1)t_2 - 1 = 2 \cdot t_2 + 1 \cdot [(q - 1)t_2 - 1] \in H_2,$$
    we conclude that $H$ is a telescopic semigroup. From the genus formula, we have
    $$\begin{array}{rl}
    g(H) & = \displaystyle{\frac{1}{2}\left(1 - a_1 + \left(\frac{d_1}{d_2} - 1\right)a_2 + \left(\frac{d_2}{d_3} - 1\right)a_3\right)} \\
    & = \displaystyle{\frac{1}{2}\left(1 - (q + 1)t_2 + \left(t_2 - 1\right)[(q^2 - 1)t_2 - (q + 1)] + q\cdot [(q + 1)t_2 - 1]\right)} \\
    & = \displaystyle{1 + \frac{[(q - 1)t_2 - 1](q + 1)t_2}{2}} \\
    & = g(\mathcal{C}_2).
    \end{array}$$
    Therefore, we conclude that $H=H(P_{\alpha}).$
\end{proof}

\section{The curve \texorpdfstring{$\mathcal{C}_3$}{C3}}\label{sec:C3}
As we mentioned before, we consider a constant field extension and analyze some places of $\overline{\mathbb{F}}_{q^n}(\mathcal{C}_3)$. Although the desingularization of the curve was already studied in \cite{guardieiro}, we include this discussion for the sake of completeness.

To study the branches centered at the singular points of $\mathcal{C}_3$, we proceed as in \cite[Section 7.4]{fischer}: the lower convex hull of the Newton polygon of the defining polynomial of $\mathcal{C}_3$ is the line passing through the points $((q^2 - 1)t_3, 0)$ and $(0, (q - 1)t_3)$, which corresponds to the polynomial $x^{(q^2 - 1)t_3} + y^{(q - 1)t_3}$. We can now determine the branches centered at the point $(0 : 0 : 1)$ by analyzing the Puiseux series. 

Since $\textrm{char }\overline{\mathbb{F}}_{q^n} = p > 0$, we need to use \cite[Theorem 5.14]{abhyankar} to ensure that such series exist. At the end of this analysis, we find $(q - 1)t_3$ branches centered at $(0 : 0 : 1)$, whose expansions in terms of a local parameter $s$ are
$$\begin{cases}
    X & = s \\
    Y & = \lambda s^{q + 1} + \cdots \\
    Z & = 1
\end{cases}$$
with $\lambda$ a $(q - 1)t_3$-th root of $-1$. Each branch is associated to a place of $\overline{\mathbb{F}}_{q^n}(\mathcal{C}_3)$, which we  denote by $P_Z^1, \ldots, P_Z^{(q - 1)t_3}$. By defining the coordinate functions $x = \frac{X}{Z}$ and $y = \frac{Y}{Z}$, we then have, for any $i = 1, \ldots, (q - 1)t_3$,
$$v_{Z_{\infty}^i}(x) = 1 \qquad \textrm{and} \qquad v_{Z_{\infty}^i}(y) = q + 1.$$

Observe that $\theta: (X : Y : Z) \mapsto (Z : X : Y)$ is an automorphism of $\mathcal{C}_3$, which maps $(0 : 0 : 1)$ into $(1 : 0 : 0)$. Therefore, we also find $(q - 1)t_3$ branches centered at $(1 : 0 : 0)$, which are associated to the places $P_X^1, \ldots, P_X^{(q - 1)t_3}$ of $\overline{\mathbb{F}}_{q^n}(\mathcal{C}_3)$, and whose expansions in terms of the local parameter $s$ are
$$\begin{cases}
    X & = 1 \\
    Y & = s \\
    Z & = \lambda s^{q + 1} + \cdots 
\end{cases}$$
with $\lambda$ a $(q - 1)t_3$-th root of $-1$, and this gives
$$v_{P_X^i}(x) = - (q + 1) \qquad \textrm{and} \qquad v_{P_X^i}(y) = -q, \qquad i = 1, \ldots, (q - 1)t_3.$$

Finally, we apply $\theta$ once more to find $(q - 1)t_3$ branches centered at $(0 : 1 : 0)$, which are associated to the places $P_Y^1, \ldots, P_Y^{(q - 1)t_3}$, whose expansions in terms of the local parameter $s$ are
$$\begin{cases}
    X & = \lambda s^{q + 1} + \cdots \\
    Y & = 1 \\
    Z & = s 
\end{cases}$$
with $\lambda$ a $(q - 1)t_3$-th root of $-1$, and 
$$v_{P_Y^i}(x) = q \qquad \textrm{and} \qquad v_{P_Y^i}(y) = -1, \qquad i = 1, \ldots, (q - 1)t_3.$$

We already noted in Proposition \ref{non-singular-points} that any other place $P_{(a,b)}$ of $\overline{\mathbb{F}}_{q^n}(\mathcal{C}_3)$ is the common zero of $x - a$ and $y - b$, so the supports of $\divi(x)$ and $\divi(y)$ are contained in
$$\mathcal{O} \coloneqq \{P_X^1, \ldots, P_X^{(q - 1)t_3}, P_Y^1, \ldots, P_Y^{(q - 1)t_3}, P_Z^1, \ldots, P_Z^{(q - 1)t_3}\}.$$
The information on the valuations of such functions yields the following result.

\begin{proposition}\label{divisors-x-y-2}
    The principal divisors of the functions $x$ and $y$ are given by
    $$\divi(x) = \sum_{i = 1}^{(q - 1)t_3}P_Z^i + q\sum_{i = 1}^{(q - 1)t_3}P_Y^i - (q + 1)\sum_{i = 1}^{(q - 1)t_3}P_X^i$$
    and
    $$\divi(y) = (q + 1)\sum_{i = 1}^{(q - 1)t_3}P_Z^i - \sum_{i = 1}^{(q - 1)t_3}P_Y^i - q\sum_{i = 1}^{(q - 1)t_3}P_X^i$$
\end{proposition}

Analogously to the previous curve, we consider some automorphisms of $\mathcal{C}_3$ to obtain the Weierstrass semigroups of a collection of places by studying only one of them.

\begin{proposition}\label{singular-points}
    The set $\mathcal{O}$ is contained in an orbit under the action of the automorphism group of $\mathcal{C}_3$. In particular, the Weierstrass semigroup at all the places of $\mathcal{O}$ is the same.
\end{proposition}

\begin{proof}
    Consider the Fermat curve of degree $(q^3 - 1)t_3$ with equation
    $$\mathcal{F}: u^{(q^3 - 1)t_3} + v^{(q^3 - 1)t_3} + 1 = 0.$$
    By defining $x = u^qv^{-(q + 1)}$ and $y = u^{-1}v^{-q}$, we can see that $\mathcal{C}_3$ is covered by $\mathcal{F}$. Moreover, the field extension $\overline{\mathbb{F}}_{q^n}(u,v) / \overline{\mathbb{F}}_{q^n}(x,y)$ is Galois, and we can follow \cite[Section 3]{butson-stewart} to conclude that $\overline{\mathbb{F}}_{q^n}(x,y)$ is fixed by the following subgroup of $\textrm{Aut}(\mathcal{F})$:
    $$G\cong \left\{\begin{pmatrix}
        \xi^{i(q^2-1)t_3}&0&0\\0&\xi^{iq(q-1)t_3}&0\\0&0&1
    \end{pmatrix}\right\}_{i=0,\ldots,q^2+q}$$
    where $\xi$ is a primitive $(q^3 - 1)t_3$-th root of unity.

    The automorphism of $\mathcal{F}$ given by $\sigma: (u, v) \mapsto (\xi^{-1}u, \xi^{-1}v)$ induces, via the covering map, on $\overline{\mathbb{F}}_{q^n}(x, y)$ the automorphism $\theta _{\xi, \xi ^{q+1}}(x, y) \mapsto (\xi x, \xi^{q+1}y)$. Indeed, since $x = \dfrac{u^{q}}{v^{q+1}}$ and $y = \dfrac{1}{u v ^{q}}$, we have 
    $$
    \frac{(\xi^{-1}u)^q}{(\xi^{-1}v)^{q+1}}
    =\xi x,
    \qquad
    \frac{1}{(\xi^{-1}u)(\xi^{-1}v)^q}
    =\xi^{q+1}y.
    $$ 
    Observe that $\theta _{\xi, \xi ^{q+1}}$ is an automorphism of $\mathcal{C}_3$.

    The set of points of $\mathcal{F}$ arising from $X_{\infty} = (1 : 0 : 0)$ is given by
    $$A \coloneqq \{(\eta : 0 : 1),\ \eta^{(q^n - 1)} = -1\}.$$
    We claim that $\sigma$ acts transitively on $A$. Indeed, let $(\eta _1 \colon 0 \colon 1), (\eta _2 \colon 0 \colon 1) \in A$. Since $(\frac{\eta _2}{\eta _1})^{q^n-1}= 1$ and $\xi$ is a primitive $(q^n-1)$-th root of unity, there exists an integer $k$ such that $\eta _2 = \xi ^{-k}\eta_1$. Therefore, 
    $$
    \sigma ^k ( \eta _1 \colon 0 \colon 1) = (\eta _2 \colon 0 \colon 1)
    $$

    Thus, the iterated images of any element of $A$ under $\sigma$ go through all the elements of $A$, and hence through all the orbits of $A/G$. Therefore, we conclude that $\theta_{\xi, \xi^{q + 1}}$ acts transitively on the set $\{P_X^1, \ldots, P_X^{(q - 1)t_3}\}$. 

    Finally, since the projectivity $\theta: (X : Y : Z) \mapsto (Z : X : Y)$ induces on $\mathcal{C}_3$ the permutation of points
    $$X_{\infty} \overset{\theta}{\longrightarrow} Y_{\infty} \overset{\theta}{\longrightarrow} Z_{\infty} \overset{\theta}{\longrightarrow} X_{\infty},$$
    we conclude that $\theta$ permutes the sets $\{P_X^1, \ldots, P_X^{(q - 1)t_3}\}$, $\{P_Y^1, \ldots, P_Y^{(q - 1)t_3}\}$ and $\{P_Z^1, \ldots, P_Z^{(q - 1)t_3}\}$ cyclically. Therefore, $\mathcal{O}$ is contained in an orbit under the action of the automorphism group of $\mathcal{C}_3$.
\end{proof}

\begin{remark}
    As we did for the curve $\mathcal{C}_2$, we can present several automorphisms of $\mathcal{C}_3$, which are the maps
    $$\theta_{\gamma,\delta}(x,y) = (\gamma x,\delta y), \qquad \gamma^{(q^2 - 1)t_3} = \gamma^{(q - 1)t_3} \delta^{(q^2 - 1)t_3} = \delta^{(q-1)t_3}$$
    and
    $$\theta(x,y) = \left(\frac{1}{y}, \frac{x}{y}\right).$$
\end{remark}

\begin{proposition}\label{divisor-f-C3}
    There exists $\zeta \in \overline{\mathbb{F}}_{q^n}^*$ such that $\zeta^{(q-1)t_3} = -1$ and
    $$
        \operatorname{div}(f) = (q^2+q+1) \sum_{i = 2}^{(q - 1)t_3} P^i_X - M P^1_X,
    $$
    where $\displaystyle{f \coloneqq \frac{y}{x(x^q - \zeta^{-1} y^{q+1})}}$ and $M \coloneqq (q^2+q+1)\big((q-1)t_3 - 1\big)$.
\end{proposition}

\begin{proof}
Throughout the proof, let $m := (q-1)t_3$ and $N := q^2+q+1$, so that $Nm = q^n-1$ and $M = (m-1)N$. With this notation, the defining equation of $\mathcal{C}_3$ reads
$$x^{(q+1)m} + x^m y^{(q+1)m} + y^m = 0.$$

Consider the function $w := y^{q+1}x^{-q}$. A direct computation with Proposition \ref{divisors-x-y-2} gives
$$\mathrm{div}(w) = N\sum_{i=1}^{m} P_Z^i - N\sum_{i=1}^{m} P_Y^i
\quad \text{and} \quad
\mathrm{div}\big(y\,x^{-(q+1)}\big) = N\sum_{i=1}^{m} P_X^i - N\sum_{i=1}^{m} P_Y^i.$$
By dividing the defining equation of $\mathcal{C}_3$ by $x^{(q+1)m}$ we obtain $w^m + 1 = -y^m x^{-(q+1)m}$, and therefore
$$\mathrm{div}(w^m+1) = mN\sum_{i=1}^{m} P_X^i - mN\sum_{i=1}^{m} P_Y^i.$$

Since $p \nmid m$, the polynomial $T^m+1$ has $m$ distinct roots in $\overline{\mathbb{F}}_{q^n}$, and $w^m+1$ is the product of the factors $w-\overline{\zeta}$ with $\overline{\zeta}^m = -1$. Each factor $w-\overline{\zeta}$ has the same poles as $w$, so its pole divisor is $N\sum_{i=1}^{m} P_Y^i$, which has degree $mN$. Hence its zero divisor also has degree $mN$. The zeros of $w-\overline{\zeta}$ are zeros of $w^m+1$, so they lie in $\{P_X^1, \dots, P_X^m\}$. Moreover, two distinct factors have no common zero, since their difference is a nonzero constant. Comparing with $\mathrm{div}(w^m+1)$, we see that each $P_X^i$ is a zero of exactly one factor, with multiplicity $mN$. As the zero divisor of each factor has degree $mN$, each factor $w-\overline{\zeta}$ vanishes at exactly one place $P_X^i$.

Let $\zeta$ be such that $w - \zeta$ is the factor associated to $P_X^1$. In particular,
$$\mathrm{div}(w-\zeta) = mN\,P_X^1 - N\sum_{i=1}^{m} P_Y^i.$$
Finally, a direct computation shows that
$$f = \frac{y}{x(x^q - \zeta^{-1}y^{q+1})} = -\zeta\, y\, x^{-(q+1)}\,(w-\zeta)^{-1},$$
and therefore
$$\mathrm{div}(f) = N\sum_{i=1}^{m} P_X^i - mN\,P_X^1 = N\sum_{i=2}^{m} P_X^i - (m-1)N\,P_X^1 = (q^2+q+1)\sum_{i=2}^{(q-1)t_3} P_X^i - M P_X^1.$$
\end{proof}

\begin{theorem}\label{semigroup3}
Let $m\coloneqq(q-1)t_3$, $N\coloneqq q^2+q+1$ and $M\coloneqq N(m-1)=(q^2+q+1)\big((q-1)t_3-1\big)$. For every place $P\in\mathcal{O}$, we have
$$H(P)=\big\langle M+qn_1+(q+1)n_2 \;:\; n_1,n_2\in\NN,\ n_1\le qn_2,\ qn_1+(q+1)n_2\le q^2+q+1\big\rangle.$$
\end{theorem}
 
\begin{proof}
By Proposition \ref{singular-points}, all the places of $\mathcal{O}$ have the same Weierstrass semigroup, so it is enough to consider $P=P_X^1$. Denote by $H$ the semigroup on the right-hand side, and let $f$ be as in Proposition \ref{divisor-f-C3}.
 
Let $n_1,n_2\in\NN$ with $n_1\le qn_2$ and $qn_1+(q+1)n_2\le N$, and consider the function $z\coloneqq x^{n_2}y^{n_1}f$. By Proposition \ref{divisors-x-y-2},
$$\operatorname{div}(x^{n_2}y^{n_1})=\big(n_2+(q+1)n_1\big)\sum_{i=1}^{m}P_Z^i+\big(qn_2-n_1\big)\sum_{i=1}^{m}P_Y^i-\big((q+1)n_2+qn_1\big)\sum_{i=1}^{m}P_X^i,$$
and, by Proposition \ref{divisor-f-C3}, $\operatorname{div}(f)=N\sum_{i=2}^{m}P_X^i-MP_X^1$. Therefore
$$\begin{aligned}
\operatorname{div}(z)={}&\big(n_2+(q+1)n_1\big)\sum_{i=1}^{m}P_Z^i+\big(qn_2-n_1\big)\sum_{i=1}^{m}P_Y^i+\big(N-qn_1-(q+1)n_2\big)\sum_{i=2}^{m}P_X^i\\
&-\big(M+qn_1+(q+1)n_2\big)P_X^1.
\end{aligned}$$
The first three coefficients are nonnegative: the first one because $n_1,n_2\ge0$, the second one because $n_1\le qn_2$, and the third one because $qn_1+(q+1)n_2\le N$. Hence $P_X^1$ is the only pole of $z$, and its pole order is $M+qn_1+(q+1)n_2$. This shows that $H\subseteq H(P_X^1)$.
 
Now let $\Lambda_q$ be the semigroup of Lemma~\ref{lem:Lambda_q}, and set $L\coloneqq N$. By Lemma~\ref{lem:Lambda_q}(a),
$$T\coloneqq\Lambda_q\cap[0,N]=\{qn_1+(q+1)n_2 \;:\; n_1,n_2\in\NN,\ n_1\le qn_2,\ qn_1+(q+1)n_2\le N\},$$
and so $H=\langle (m-1)N+s \;:\; s\in T\rangle$. By Lemma \ref{lem:Lambda_q}(c), $c(\Lambda_q)+\mu(\Lambda_q)-1=N$ and $g(\Lambda_q)=\frac{q(q+1)}{2}$. Applying Proposition~\ref{lem:genus} with $L=N$ and $a=m$, and using that $N-q(q+1)=1$, we obtain
$$\begin{aligned}
g(H)&=\frac{m(m-1)N}{2}-(m-1)+m\cdot\frac{q(q+1)}{2}\\
&=\frac{Nm^2-m\big(N-q(q+1)\big)-2m+2}{2}\\&=\frac{Nm^2-3m+2}{2}.
\end{aligned}$$
On the other hand, $(q^3-1)(q-1)t_3^2=(q^2+q+1)(q-1)^2t_3^2=Nm^2$, so Proposition \ref{curves-places} gives
$$g(\Ccal_3)=\frac{(q^3-1)(q-1)t_3^2-3(q-1)t_3+2}{2}=\frac{Nm^2-3m+2}{2}=g(H).$$
Since $H\subseteq H(P_X^1)$ and both semigroups have genus $g(\Ccal_3)$ we have $H=H(P_X^1)$.
\end{proof}

\section{AI Disclosure}
This paper is the outcome of a project that started in 2024. Section \ref{sec:C2} was written back then, with no use of AI at all. For Section \ref{sec:C3} we used Claude to provide explicit functions for the cases $(q,n) \in \{(2,3),(2,6),(3,3),(3,6)\}$. From these four examples we guessed the general form of the function $f$ and then proved Proposition \ref{divisor-f-C3} ourselves.

All mathematical content, arguments, and results are the authors’ own; the authors have independently verified the final manuscript and assume full responsibility for its content.

\section{Acknowledgments}
Tiago Aprigio was supported by São Paulo Research Foundation (FAPESP) – Brazil, grant 2023/07508-1. João Paulo Guardieiro was partially supported by Coordenação de Aperfeiçoamento de Pessoal de Nível Superior – Brasil (CAPES) – Finance Code 001. Leonardo Oliveira was supported by Coordenação de Aperfeiçoamento de Pessoal de Nível Superior – Brasil (CAPES) – Finance Code 001.

We thank Lara Vicino for discussions with the second author on some of the results
of this paper.

\bibliographystyle{alpha}
\bibliography{refs}

@article{abq2019,
  author  = {Abd{\'o}n, M. and Borges, H. and Quoos, L.},
  title   = {{W}eierstrass points on {K}ummer extensions},
  journal = {Adv. Geom.},
  volume  = {19},
  number  = {3},
  year    = {2019},
  pages   = {323--333},
  doi     = {10.1515/advgeom-2018-0021}
}

@book{abhyankar,
  author    = {Abhyankar, S. S.},
  title     = {Lectures on expansion techniques in algebraic geometry},
  series    = {Tata Institute of Fundamental Research Lectures on Mathematics and Physics},
  volume    = {57},
  note      = {Notes by Balwant Singh},
  publisher = {Tata Institute of Fundamental Research},
  address   = {Bombay},
  year      = {1977}
}

@article{aprigio,
  author  = {Aprigio, T. and Guardieiro, J. P.},
  title   = {Minimal value set binomials and {F}robenius nonclassical curves},
  journal = {Finite Fields Appl.},
  volume  = {117},
  year    = {2027},
  pages   = {Paper No. 102875}
}

@article{vicino1,
  author  = {Beelen, P. and Montanucci, M. and Vicino, L.},
  title   = {{W}eierstrass semigroups and automorphism group of a maximal curve with the third largest genus},
  journal = {Finite Fields Appl.},
  volume  = {92},
  year    = {2023},
  pages   = {Paper No. 102300}
}

@article{vicino2,
  author  = {Beelen, P. and Montanucci, M. and Vicino, L.},
  title   = {{W}eierstrass semigroups and automorphism group of a maximal function field with the third largest possible genus, $q\equiv 0 \pmod{3}$},
  journal = {Finite Fields Appl.},
  volume  = {110},
  year    = {2026},
  pages   = {Paper No. 102729}
}

@article{vicino3,
  author  = {Beelen, P. and Montanucci, M. and Vicino, L.},
  title   = {{W}eierstrass semigroups and automorphism group of a maximal function field with the third largest possible genus, $q\equiv 1 \pmod{3}$},
  journal = {Finite Fields Appl.},
  volume  = {109},
  year    = {2026},
  pages   = {Paper No. 102701}
}

@article{butson-stewart,
  author  = {Butson, A. T. and Stewart, B. M.},
  title   = {Systems of linear congruences},
  journal = {Canad. J. Math.},
  volume  = {7},
  year    = {1955},
  pages   = {358--368}
}

@book{fischer,
  author    = {Fischer, G.},
  title     = {Plane algebraic curves},
  series    = {Student Mathematical Library},
  volume    = {15},
  publisher = {American Mathematical Society},
  address   = {Providence, RI},
  year      = {2001}
}

@article{goppa,
  author  = {Goppa, V. D.},
  title   = {Codes on algebraic curves},
  journal = {Soviet Math. Dokl.},
  volume  = {24},
  year    = {1981},
  pages   = {170--172}
}

@phdthesis{guardieiro,
  author  = {Guardieiro, J. P.},
  title   = {Counting rational points with {S}t{\"o}hr--{V}oloch theory and {T}ate--{S}hafarevich results},
  school  = {Universidade de S{\~a}o Paulo},
  address = {S{\~a}o Carlos},
  year    = {2025},
  note    = {\url{https://doi.org/10.11606/T.55.2025.tde-05082025-170919}}
}

@article{hamming,
  author  = {Hamming, R. W.},
  title   = {Error detecting and error correcting codes},
  journal = {Bell System Tech. J.},
  volume  = {29},
  number  = {2},
  year    = {1950},
  pages   = {147--160}
}

@article{hasse,
  author  = {Hasse, H.},
  title   = {{\"U}ber den algebraischen {F}unktionenk{\"o}rper der {F}ermatschen {G}leichung},
  journal = {Acta Sci. Math. (Szeged)},
  volume  = {13},
  year    = {1950},
  pages   = {195--207}
}

@article{hlp1998,
  author  = {Kirfel, C. and Pellikaan, R.},
  title   = {The minimum distance of codes in an array coming from telescopic semigroups},
  journal = {IEEE Trans. Inform. Theory},
  volume  = {41},
  number  = {6},
  year    = {1995},
  pages   = {1720--1732}
}

@book{numerical,
  author    = {Rosales, J. C. and Garc{\'\i}a-S{\'a}nchez, P. A.},
  title     = {Numerical semigroups},
  series    = {Developments in Mathematics},
  volume    = {20},
  publisher = {Springer},
  address   = {New York},
  year      = {2009}
}

@book{stichtenoth,
  author    = {Stichtenoth, H.},
  title     = {Algebraic function fields and codes},
  edition   = {Second},
  series    = {Graduate Texts in Mathematics},
  volume    = {254},
  publisher = {Springer-Verlag},
  address   = {Berlin},
  year      = {2009}
}

@article{stohr-voloch,
  author  = {St{\"o}hr, K.-O. and Voloch, J. F.},
  title   = {{W}eierstrass points and curves over finite fields},
  journal = {Proc. London Math. Soc. (3)},
  volume  = {52},
  number  = {1},
  year    = {1986},
  pages   = {1--19}
}

@article{sylvester1884,
  author  = {Sylvester, J. J.},
  title   = {Mathematical questions with their solutions},
  journal = {Educational Times},
  volume  = {41},
  year    = {1884},
  pages   = {21}
}

@article{tsfasman,
  author  = {Tsfasman, M. A. and Vl{\u{a}}du{\c{t}}, S. G. and Zink, Th.},
  title   = {Modular curves, {S}himura curves, and {G}oppa codes, better than {V}arshamov--{G}ilbert bound},
  journal = {Math. Nachr.},
  volume  = {109},
  year    = {1982},
  pages   = {21--28}
}

@article{yang-hu,
  author  = {Yang, S. and Hu, C.},
  title   = {{W}eierstrass semigroups from {K}ummer extensions},
  journal = {Finite Fields Appl.},
  volume  = {45},
  year    = {2017},
  pages   = {264--284}
}

@misc{zhang-zhao,
  author       = {Zhang, H. and Zhao, C.-A.},
  title        = {{W}eierstrass semigroups at totally ramified places of degree one on {K}ummer extensions},
  howpublished = {preprint, \texttt{arXiv:2605.14583 [math.AG]}},
  year         = {2026}
}

@article{castellanos2024,
  author  = {Castellanos, A. S. and Mendoza, E. A. R. and Quoos, L.},
  title   = {{W}eierstrass semigroups, pure gaps and codes on function fields},
  journal = {Des. Codes Cryptogr.},
  volume  = {92},
  number  = {5},
  pages   = {1219--1242},
  year    = {2024}
}

\end{document}